\documentclass[11pt,a4paper]{article}

\usepackage[T1]{fontenc}
\usepackage{lmodern}
\usepackage[a4paper,margin=2.8cm]{geometry}
\usepackage{amsmath,amssymb,amsthm,mathtools}
\usepackage{enumitem}
\usepackage{microtype}
\usepackage[hidelinks]{hyperref}
\hypersetup{
  pdftitle={Comparing the Automorphism Group of the Measure Algebra with Some Groups Related to the Infinite Permutation Group of the Natural Numbers},
  pdfauthor={Pietro Ursino},
  pdfsubject={Automorphism groups of measure algebras and quotient Boolean algebras},
  pdfkeywords={measure algebra, Boolean algebra, automorphism group, locally moving groups}
}

\allowdisplaybreaks
\newtheorem{theorem}{Theorem}[section]
\newtheorem{lemma}[theorem]{Lemma}
\newtheorem{corollary}[theorem]{Corollary}
\newtheorem{proposition}[theorem]{Proposition}
\theoremstyle{definition}
\newtheorem{definition}[theorem]{Definition}
\theoremstyle{remark}

\newcommand{\N}{\mathbb N}
\newcommand{\B}{\mathbb B}
\newcommand{\M}{\mathcal M}
\newcommand{\id}{\mathrm{Id}}
\newcommand{\var}{\operatorname{var}}
\newcommand{\fix}{\operatorname{fix}}
\newcommand{\Sp}{\operatorname{Sp}}
\newcommand{\Aut}{\operatorname{Aut}}
\newcommand{\MPAut}{\operatorname{MPAut}}
\newcommand{\Z}{\operatorname{Z}}
\newcommand{\Pow}{\mathcal P}
\newcommand{\fin}{\mathrm{fin}}
\newcommand{\set}[1]{\left\{#1\right\}}

\newcommand{\meet}{\wedge}
\newcommand{\join}{\vee}

\title{Comparing the Automorphism Group of the Measure Algebra with Some Groups Related to the Infinite Permutation Group of the Natural Numbers}
\author{Pietro Ursino\\
\small Department of Mathematics and Informatics\\
  \small Università degli Studi di Catania\\
  \small Viale Andrea Doria 6, 95125 Catania, Italy}

\date{}

\begin{document}
\maketitle

\begin{abstract}
Let $\M$ be the standard Lebesgue measure algebra, let $\Aut(\M)$ be its Boolean-algebra automorphism group, and let $\MPAut(\M)$ be the subgroup of measure-preserving automorphisms. We prove that neither of these groups is isomorphic to either of the groups $T_\omega^*$ and $S_\omega^*$ associated with the trivial automorphisms of $\Pow(\N)/\fin$. The proof defines the support preorder by a formula in the pure group language. Any hypothetical group isomorphism would therefore induce an order embedding of $\Pow(\N)/\fin$ into $\M$, contradicting the countable chain condition of the measure algebra.
\end{abstract}



\section{Introduction}

This paper compares the automorphism group of the standard Lebesgue measure algebra with groups associated with the quotient Boolean algebra $\Pow(\N)/\fin$. The measure algebra $\M$ is the Boolean algebra of Lebesgue measurable subsets of $[0,1]$ modulo null sets. We write $\Aut(\M)$ for its full Boolean-algebra automorphism group and $\MPAut(\M)$ for the subgroup of automorphisms that preserve Lebesgue measure; background on automorphism groups of measure algebras can be found in Fremlin\cite{fremlin}.

A \emph{near permutation} of $\N$ is a bijection $p:D\to R$ between cofinite subsets $D,R\subseteq\N$. It induces an automorphism of $\Pow(\N)/\fin$ by
\[
  [A]\longmapsto [p[A\cap D]],
\]
where $[A]$ denotes the equivalence class of $A$ modulo finite sets. We use the standard notation $T_\omega^*$ for the group of trivial automorphisms induced by near permutations, and $S_\omega^*$ for the subgroup induced by permutations of all of $\N$ (permutations that differ only on a finite set induce the same automorphism).

The main result is the following.

\begin{theorem}\label{thm:main-intro}
Neither $\Aut(\M)$ nor $\MPAut(\M)$ is isomorphic, as an abstract group, to $T_\omega^*$ or to $S_\omega^*$.
\end{theorem}

The comparison is related to the set-theoretic problem of embedding the measure algebra into $\Pow(\N)/\fin$, which has independent axiomatic subtleties; see Dow and Hart\cite{dowhart}. The proof reconstructs the inclusion relation between supports from the pure group structure. This is in the spirit of the reconstruction methods developed by Rubin and others; see, for example, \cite{rubin1989,rubin1996}. The resulting order embedding is incompatible with the order structure of the measure algebra: $\Pow(\N)/\fin$ has a strictly increasing chain of order type $\omega_1$, while $\M$ satisfies the countable chain condition.

\section{Supports and elementary facts}

We use the Boolean-algebra operations $\meet$, $\join$, and complement $\neg a$. The order is the usual Boolean order. Unless otherwise indicated, $\B$ denotes either $\M$ or $\Pow(\N)/\fin$, and $G$ denotes one of $\Aut(\M)$, $\MPAut(\M)$, $T_\omega^*$, or $S_\omega^*$, acting on the corresponding algebra.

\begin{definition}\label{def:support}
For an automorphism $f$ of a Boolean algebra $\B$, whenever the indicated suprema exist, define
\[
 \var(f)=\bigvee\set{a\in\B:f(a)\meet a=0}
\]
and
\[
 \fix(f)=\bigvee\set{a\in\B:(\forall b\leq a)\ f(b)=b}.
\]
For a subgroup $G\leq\Aut(\B)$ and $a\in\B$, put
\[
 \Sp_G(a)=\set{g\in G:\var(g)\leq a},\qquad
 \Sp_G^*(a)=\set{g\in G:\var(g)=a}.
\]
For $X\subseteq G$, define
\[
 \Z(X)=\set{g\in G:(\forall x\in X)\ gx=xg},
 \qquad \Z(f)=\Z(\{f\}).
\]
We use the commutator and conjugation conventions
\[
 [f,g]=fgf^{-1}g^{-1},\qquad f^h=hfh^{-1},\qquad [h,f,g]=[[h,f],g].
\]
\end{definition}

The indicated suprema exist automatically when $\B$ is complete. For complete Boolean algebras, $\neg\var(f)=\fix(f)$. The same identity holds in $\Pow(\N)/\fin$ for the automorphisms considered here, as the following lemma shows.

\begin{lemma}\label{lem:near-permutation-support}
Let $p:D\to R$ be a near permutation, where $D,R\subseteq\N$ are cofinite, and let $f_p$ be the induced automorphism of $\Pow(\N)/\fin$. Set
\[
 E_p=\set{n\in D:p(n)\ne n},\qquad F_p=\set{n\in D:p(n)=n}.
\]
Then the supremum defining $\var(f_p)$ exists in $\Pow(\N)/\fin$ and
\[
 \var(f_p)=[E_p],\qquad \fix(f_p)=[F_p].
\]
In particular, these classes do not depend on the chosen representative of the trivial automorphism.
\end{lemma}

\begin{proof}
Let
\[
 \mathcal D_f=\set{[A]\in\Pow(\N)/\fin:p[A\cap D]\cap A\text{ is finite}}.
\]
If $[A]\in\mathcal D_f$, then $A\cap F_p$ is finite, since each point of $A\cap F_p$ belongs to $p[A\cap D]\cap A$. As $D$ is cofinite, $[E_p]=[\N\setminus F_p]$, so $[A]\leq[E_p]$. Thus $[E_p]$ is an upper bound for $\mathcal D_f$.

To prove it is the least upper bound, suppose $[C]$ is an upper bound and $[E_p]\nleq[C]$. Then $E_p\setminus C$ is infinite. Form an undirected graph on $\N$ by joining $n$ to $p(n)$ whenever $n\in D$ and $p(n)\ne n$. Each vertex has degree at most two, because $p$ is a partial bijection. Every infinite subset of a graph of finite degree has an infinite independent subset: recursively choose a vertex and discard it and its finitely many neighbours. Hence there is an infinite independent set $B\subseteq E_p\setminus C$. Since $p[B]\cap B=\varnothing$, we have $[B]\in\mathcal D_f$, but $[B]\nleq[C]$, a contradiction. Therefore $[E_p]$ is the supremum of $\mathcal D_f$.

If $[A]\leq[F_p]$, then all but finitely many points of $A$ are fixed by $p$, and so $f_p([B])=[B]$ for every $[B]\leq[A]$. Thus $[A]\leq\fix(f_p)$. Conversely, if $[A]\nleq[F_p]$, then $A\cap E_p$ is infinite. The same graph argument gives an infinite independent set $B\subseteq A\cap E_p$. Then $[B]\leq[A]$ but $f_p([B])\ne[B]$, so $[A]\nleq\fix(f_p)$. This proves $\fix(f_p)=[F_p]$ and completes the proof.
\end{proof}

We will use the following standard support identities:
\begin{align}
 \var(u^h)&=h(\var(u)),\label{eq:conj-support}\\
 \var(uv)&\leq\var(u)\join\var(v),\label{eq:prod-support}\\
 \var([u,v])&\leq\var(u)\join v(\var(u)).\label{eq:comm-support}
\end{align}
Automorphisms with disjoint variations commute. These identities hold in $\M$ and, by Lemma~\ref{lem:near-permutation-support}, in $\Pow(\N)/\fin$ for the groups under consideration.

\begin{proposition}\label{prop:support-properties}
For each of the four groups $G$ under consideration the following statements hold.
\begin{enumerate}[label=(\alph*),leftmargin=2em]
\item If $0\ne a\leq\bigwedge_{i=1}^m\var(g_i)$, then there exists $0\ne b\leq a$ such that
\[
 b\meet g_i(b)=0\qquad(1\leq i\leq m).
\]
\item For every $0\ne a\in\B$ and every integer $n\geq1$, there are $h\in G$ and $0\ne c\leq a$ such that $\var(h)=a$ and
\[
 c,h(c),\ldots,h^n(c)
\]
are pairwise disjoint. In particular, $h^n\ne\id$.
\end{enumerate}
\end{proposition}

\begin{proof}
For $\B=\M$, part (a) follows by a finite iteration. If $0\ne d\leq\var(g_i)$, the definition of the supremum gives an element $e$ with $e\meet g_i(e)=0$ and $d\meet e\ne0$. Replacing $d$ by $d\meet e$ gives a nonzero subelement disjoint from its image under $g_i$. Repeating this for $i=1,\ldots,m$ preserves all disjointness properties already obtained.

For $\B=\Pow(\N)/\fin$, choose an infinite representative $A$ of $a$ and near-permutation representatives $p_i:D_i\to R_i$ of the $g_i$. Since $a\leq\var(g_i)$, after deleting finitely many points from $A$, every remaining point belongs to every $D_i$ and is moved by every $p_i$. Put an edge between $x$ and $p_i(x)$ whenever $p_i(x)\ne x$. The union graph has finite degree at most $2m$. It therefore has an infinite independent subset $C$ of the remaining part of $A$. Then $p_i[C]\cap C=\varnothing$ for every $i$, and $b=[C]$ proves (a).

For part (b) in $\Pow(\N)/\fin$, enumerate an infinite representative of $a$ as $A=\{a_k:k\in\mathbb Z\}$. Let $h$ be induced by the permutation that sends $a_k$ to $a_{k+1}$ and fixes $\N\setminus A$. By Lemma~\ref{lem:near-permutation-support}, $\var(h)=a$. The class of $C=\{a_{(n+1)k}:k\in\mathbb Z\}$ has pairwise disjoint translates $C,h(C),\ldots,h^n(C)$.

For $\B=\M$, the interval below any nonzero $a$ is, with its normalized measure, isomorphic to the Lebesgue measure algebra of the circle. On that interval take an irrational rotation and extend it by the identity on the complement of $a$. The resulting measure-preserving automorphism $h$ has variation $a$. The first $n+1$ iterates of a point under an irrational rotation are distinct; hence a sufficiently short arc represents a nonzero $c\leq a$ whose first $n+1$ translates are pairwise disjoint. This same automorphism belongs to both $\Aut(\M)$ and $\MPAut(\M)$.
\end{proof}

\section{A group-theoretic definition of the support order}

We now define a preorder using only the group operations. Thus, although the definition refers to the group $G$, it is preserved by every abstract group isomorphism between groups for which the definition is used.

\begin{definition}\label{def:phi}
For $f,f'\in G$, let
\begin{align*}
 \phi_1(f,f')\iff{}&\ \forall g\in G\Bigl([g,f]\ne\id\ \Longrightarrow\ \exists f_1,f_2\in\Z(f')\text{ such that}\\[-0.4em]
 &\hspace{6em}[g,f_1,f_2]\ne\id\quad\text{and}\quad [[g,f_1,f_2],f']=\id\Bigr).
\end{align*}
Put
\[
 D_1(f)=\set{f'\in G:\phi_1(f,f')},\qquad
 V(f)=\Z\bigl(\set{(f')^{12}:f'\in D_1(f)}\bigr).
\]
Finally, define
\[
 \phi_{\leq}(f,g)\iff V(f)\subseteq V(g),\qquad
 \phi_{=}(f,g)\iff \phi_{\leq}(f,g)\text{ and }\phi_{\leq}(g,f).
\]
\end{definition}

The power $12$ is chosen so that the support of the twelfth power is contained in the support of each of the first, second, third, and fourth powers. More generally, for any automorphism $u$ and any positive integer $i$ dividing $12$,
\[
 \var(u^{12})\leq\var(u^i).
\]
Indeed, every Boolean-algebra interval fixed pointwise by $u^i$ is also fixed pointwise by $u^{12}=(u^i)^{12/i}$.

\begin{lemma}\label{lem:phi1-characterization}
For each of the four groups $G$:
\begin{enumerate}[label=(\alph*),leftmargin=2em]
\item If $\var(f)\meet\var(f')=0$, then $\phi_1(f,f')$.
\item If $\var(f)\meet\var((f')^{12})\ne0$, then $\neg\phi_1(f,f')$.
\end{enumerate}
\end{lemma}

The proof of Lemma~\ref{lem:phi1-characterization}, together with the reconstruction of the support order, is given in Appendix~\ref{app:proof}.

\begin{theorem}\label{thm:support-order}
For $G=T_\omega^*$ or $G=S_\omega^*$, the predicate $\phi_{\leq}(f,g)$ holds if and only if $\var(f)\leq\var(g)$ in $\Pow(\N)/\fin$. For $G=\Aut(\M)$ or $G=\MPAut(\M)$, it holds if and only if $\var(f)\leq\var(g)$ in $\M$.
\end{theorem}

\begin{corollary}\label{cor:support-equality}
For any of the four groups,
\[
 \phi_{=}(f,g)\quad\Longleftrightarrow\quad\var(f)=\var(g).
\]
\end{corollary}

\section{Non-isomorphism results}

We first record the order-theoretic fact needed in the proof.

\begin{lemma}\label{lem:omega1-chain}
The Boolean algebra $\Pow(\N)/\fin$ contains a strictly increasing chain of order type $\omega_1$.
\end{lemma}

\begin{proof}
We construct by transfinite recursion infinite sets $A_\alpha\subseteq\N$ for $\alpha<\omega_1$ such that, whenever $\alpha<\beta<\omega_1$,
\[
 A_\beta\subseteq^* A_\alpha\qquad\text{and}\qquad A_\alpha\setminus A_\beta\text{ is infinite},
\]
where $\subseteq^*$ denotes inclusion modulo a finite set. Start with $A_0=\N$. At a successor stage, choose an infinite subset $A_{\alpha+1}\subseteq A_\alpha$ with infinite complement in $A_\alpha$.

Let $\lambda<\omega_1$ be a nonzero limit ordinal, and choose a strictly increasing sequence $\gamma_0<\gamma_1<\cdots$ cofinal in $\lambda$. For each $n$, the set
\[
 \left(\bigcap_{j\leq n} A_{\gamma_j}\right)\setminus A_{\gamma_{n+1}}
\]
is infinite: the finite intersection differs from $A_{\gamma_n}$ by at most finitely many points, and $A_{\gamma_n}\setminus A_{\gamma_{n+1}}$ is infinite. Choose distinct points
\[
 x_n\in\left(\bigcap_{j\leq n} A_{\gamma_j}\right)\setminus A_{\gamma_{n+1}}
\]
and set $A_\lambda=\{x_n:n\in\N\}$. For every fixed $m$, all $x_n$ with $n\geq m$ lie in $A_{\gamma_m}$, so $A_\lambda\subseteq^*A_{\gamma_m}$. Moreover, $A_{\gamma_m}\setminus A_\lambda$ is infinite, since the infinite set $A_{\gamma_m}\setminus A_{\gamma_{m+1}}$ contains at most finitely many points of $A_\lambda$. Cofinality gives the required two properties for every $\alpha<\lambda$. The construction is complete. Taking complements gives a strictly increasing $\omega_1$-chain in $\Pow(\N)/\fin$.
\end{proof}

\begin{proof}[Proof of Theorem~\ref{thm:main-intro}]
Let $H$ be either $T_\omega^*$ or $S_\omega^*$, and let $G$ be either $\Aut(\M)$ or $\MPAut(\M)$. Suppose for a contradiction that $\Phi:H\to G$ is a group isomorphism.

For every $a\in\Pow(\N)/\fin$, choose $f_a\in H$ with $\var(f_a)=a$; such an element exists by Proposition~\ref{prop:support-properties}(b), with the case $a=0$ given by the identity. Define
\[
 \Theta(a)=\var(\Phi(f_a))\in\M.
\]
This is independent of the choice of $f_a$. Indeed, if $\var(f_a)=\var(f'_a)$, Corollary~\ref{cor:support-equality} gives $\phi_{=}(f_a,f'_a)$. Since $\phi_{=}$ is expressed in the pure group language, it is preserved by $\Phi$, and another application of Corollary~\ref{cor:support-equality} gives $\var(\Phi(f_a))=\var(\Phi(f'_a))$.

For any $a,b\in\Pow(\N)/\fin$, Theorem~\ref{thm:support-order} and preservation of the group-theoretic predicate $\phi_{\leq}$ imply
\[
 a\leq b\ \Longleftrightarrow\ \phi_{\leq}(f_a,f_b)
 \ \Longleftrightarrow\ \phi_{\leq}(\Phi(f_a),\Phi(f_b))
 \ \Longleftrightarrow\ \Theta(a)\leq\Theta(b).
\]
Thus $\Theta$ is an order embedding. By Lemma~\ref{lem:omega1-chain}, $\M$ would contain a strictly increasing chain of order type $\omega_1$. This is impossible: $\M$ satisfies the countable chain condition, whereas the successive differences of a strictly increasing $\omega_1$-chain would be an uncountable family of pairwise disjoint nonzero elements. This contradiction proves the theorem.
\end{proof}

A subgroup $G\leq\Aut(\B)$ of a complete atomless Boolean algebra is called \emph{locally moving} if, for every $0\ne a\in\B$, there is a nonidentity $g\in G$ with $\var(g)\leq a$.

\begin{theorem}\label{thm:rubin-generalization}
Let $\B$ be a complete atomless Boolean algebra with no strictly increasing chain of order type $\omega_1$, and let $G\leq\Aut(\B)$ be locally moving. Then $G$ is isomorphic to neither $T_\omega^*$ nor $S_\omega^*$.
\end{theorem}

\begin{proof}
We use Rubin's reconstruction theorem: a complete atomless Boolean algebra can be reconstructed from any locally moving subgroup of its automorphism group; in particular, an abstract isomorphism between two such locally moving groups yields an isomorphism between the underlying complete Boolean algebras. See \cite{rubin1996,mcclearyrubin}.

Put $\B_0=\Pow(\N)/\fin$ and let $\widehat{\B}_0$ be its Boolean completion. Every automorphism of $\B_0$ extends uniquely to an automorphism of $\widehat{\B}_0$. Let $\widehat T_\omega^*$ and $\widehat S_\omega^*$ be the extension groups corresponding to $T_\omega^*$ and $S_\omega^*$. They are isomorphic to their respective original groups. They are locally moving: if $0\ne c\in\widehat{\B}_0$, choose $0\ne a\in\B_0$ with $a\leq c$, using the order density of $\B_0$ in its completion. By Proposition~\ref{prop:support-properties}(b), there is a nonidentity permutation-induced automorphism with variation $a$; its extension is nontrivial and is supported below $a$, hence below $c$.

If $G$ were isomorphic to $T_\omega^*$ or $S_\omega^*$, Rubin's reconstruction theorem would imply $\B\cong\widehat{\B}_0$. But $\B_0$ has a strictly increasing $\omega_1$-chain by Lemma~\ref{lem:omega1-chain}, and its inclusion in $\widehat{\B}_0$ preserves order. The isomorphism would transfer this chain to $\B$, contradicting the hypothesis on $\B$.
\end{proof}

\appendix
\section{Proof of Theorem~\ref{thm:support-order}}\label{app:proof}

In this appendix, $G$ is any of the four groups specified in Section~2.

\begin{lemma}\label{lem:centralizer-cover}
Let $k_0,\ldots,k_n\in\mathbb Z$, $f\in\Aut(\B)$, and $a\in\B$. Suppose that $f^{k_0}(a),\ldots,f^{k_n}(a)$ are pairwise disjoint. If $h_1,\ldots,h_n\in\Z(f)$ and $0\ne b\leq a$, then
\[
 \bigvee_{i=0}^{n}f^{k_i}(b)\nleq\bigvee_{i=1}^{n}h_i(a).
\]
\end{lemma}

\begin{proof}
We argue by induction on $n$. For $n=1$, suppose that
\[
 f^{k_0}(b)\vee f^{k_1}(b)\leq h_1(a).
\]
Applying $f^{-k_0}$ and $f^{-k_1}$ shows that $b$ is below both $f^{-k_0}h_1(a)$ and $f^{-k_1}h_1(a)$. Their meet is
\[
 h_1\bigl(f^{-k_0}(a)\meet f^{-k_1}(a)\bigr)=0,
\]
because $h_1$ commutes with $f$ and the two displayed translates of $a$ are disjoint. This contradicts $b\ne0$.

Suppose the result holds for $n$ and consider $n+2$ translates, indexed by $0,\ldots,n+1$, and $n+1$ centralizers. Write
\[
 S=\bigvee_{i=1}^{n}h_i(a),\qquad H=h_{n+1}(a),
\]
and suppose that $\bigvee_{i=0}^{n+1} f^{k_i}(b)\leq S\vee H$. If $H$ is disjoint from every $f^{k_i}(b)$, then each translate is below $S$, contradicting the induction hypothesis. Otherwise choose $j$ with $H\meet f^{k_j}(b)\ne0$, and replace $b$ by the nonzero element
\[
 d=f^{-k_j}\bigl(H\meet f^{k_j}(b)\bigr)\leq b.
\]
Then $f^{k_j}(d)\leq H$.

For each $i\ne j$, restrict the current nonzero element $d$ further. The base case, applied to the exponents $k_i,k_j$ and to $h_{n+1}$, implies that $f^{k_i}(d)\vee f^{k_j}(d)\nleq H$. Since $f^{k_j}(d)\leq H$, the element $f^{k_i}(d)\setminus H$ is nonzero. Replace $d$ by
\[
 f^{-k_i}\bigl(f^{k_i}(d)\setminus H\bigr).
\]
This replacement is a nonzero subelement of the previous one, its $k_i$-translate is disjoint from $H$, and its $k_j$-translate remains below $H$. After finitely many such restrictions, we have $0\ne d\leq b$ such that $f^{k_j}(d)\leq H$ and $f^{k_i}(d)\meet H=0$ for every $i\ne j$.

The assumed cover now implies $f^{k_i}(d)\leq S$ for every $i\ne j$. Applying the induction hypothesis to the $n+1$ translates indexed by $i\ne j$ and to $h_1,\ldots,h_n$ gives a contradiction. This completes the induction.
\end{proof}

\begin{proof}[Proof of Lemma~\ref{lem:phi1-characterization}]
(a) Suppose $\var(f)\meet\var(f')=0$ and let $g\in G$ satisfy $[g,f]\ne\id$. Automorphisms with disjoint variations commute, so $\var(g)\meet\var(f)\ne0$. By Proposition~\ref{prop:support-properties}(a), choose $0\ne a\leq\var(g)\meet\var(f)$ such that $g(a)\meet a=0$. Choose $f_1\in G$ with $\var(f_1)=a$ and set $g_1=[g,f_1]$. Since $f_1^g=gf_1g^{-1}$ and $f_1^{-1}$ have disjoint variations $g(a)$ and $a$, respectively,
\[
 \var(g_1)=a\vee g(a),\qquad g_1(a)=a.
\]
In particular, $a\leq\var(g_1)$. Choose $0\ne b\leq a$ such that $g_1(b)\meet b=0$, and choose $f_2\in G$ with $\var(f_2)=b$. Then
\[
 g_2=[g_1,f_2]\ne\id,\qquad \var(g_2)=b\vee g_1(b)\leq a.
\]
Both $f_1$ and $f_2$ have variation below $\var(f)$ and therefore commute with $f'$. Also $\var(g_2)\leq\var(f)$, so $g_2$ commutes with $f'$. This supplies the witnesses required by $\phi_1(f,f')$.

(b) Suppose $a=\var(f)\meet\var((f')^{12})\ne0$. Since $12$ is divisible by $1,2,3,4$, the support identities imply
\[
 a\leq\bigwedge_{i=1}^{4}\var((f')^i).
\]
By Proposition~\ref{prop:support-properties}(a), choose $0\ne b\leq a$ such that
\[
 b\meet\bigvee_{i=1}^{4}(f')^i(b)=0.
\]
Then $b,f'(b),\ldots,(f')^4(b)$ are pairwise disjoint. Since $b\leq\var(f)$, Lemma~\ref{lem:local-moving-tools}(d) gives $g\in\Sp_G(b)$ with $[g,f]\ne\id$.

Fix arbitrary $f_1,f_2\in\Z(f')$ and let $g_2=[g,f_1,f_2]$. If $g_2=\id$, the nontriviality condition in $\phi_1$ fails. Suppose $g_2\ne\id$ and, towards a contradiction, that $g_2$ commutes with $f'$. By \eqref{eq:comm-support},
\[
 \var(g_2)\leq b\vee f_1(b)\vee f_2(b)\vee f_2f_1(b).
\]
Let $Q=\{\id,f_1,f_2,f_2f_1\}$. Since $\var(g_2)\ne0$, choose $q\in Q$ such that
\[
 c=\var(g_2)\meet q(b)\ne0,
 \qquad d=q^{-1}(c)\leq b.
\]
Every element of $Q$ commutes with $f'$. As $g_2$ also commutes with $f'$, its variation is invariant under $f'$. It follows that
\[
 (f')^i(d)\leq q^{-1}(\var(g_2))\leq\bigvee_{r\in Q}q^{-1}r(b),
 \qquad i=0,1,2,3,4.
\]
The five translates $b,f'(b),\ldots,(f')^4(b)$ are pairwise disjoint, while the four elements $q^{-1}r$, $r\in Q$, belong to $\Z(f')$. Lemma~\ref{lem:centralizer-cover}, applied to these five translates and the nonzero element $d\leq b$, says that the displayed inclusion is impossible. Thus $g_2$ cannot commute with $f'$. Since $f_1,f_2$ were arbitrary, $\phi_1(f,f')$ fails.
\end{proof}

\begin{lemma}\label{lem:local-moving-tools}
For each of the four groups $G$:
\begin{enumerate}[label=(\alph*),leftmargin=2em]
\item If $0\ne a\leq\bigwedge_{i=1}^n\var(g_i)$, there exists $0\ne b\leq a$ such that $b\meet\bigvee_{i=1}^n g_i(b)=0$.
\item For every $0\ne a\in\B$ and $n\geq1$, there is $h\in\Sp_G^*(a)$ with $h^n\ne\id$.
\item If $0\ne a\leq\var(f)\meet\var(g)$, there exists $h\in\Sp_G(a)$ such that $f^h$ does not commute with $g$.
\item If $g\in G$ and $0\ne a\leq\var(g)$, there exists $k\in\Sp_G(a)$ that does not commute with $g$.
\end{enumerate}
\end{lemma}

\begin{proof}
(a) is Proposition~\ref{prop:support-properties}(a), and (b) is Proposition~\ref{prop:support-properties}(b).

(c) If $f$ and $g$ do not commute, take $h=\id$. Otherwise, by (a), choose $0\ne b_1\leq a$ such that
\[
 f(b_1)\meet b_1=0,\qquad g(b_1)\meet b_1=0.
\]
If $b_1\meet\var(fg)\ne0$, use (a) again to choose $0\ne b\leq b_1\meet\var(fg)$ with $fg(b)\meet b=0$. Otherwise set $b=b_1$; in this case $fg$ fixes every subelement of $b$. By (b), choose $h\in\Sp_G^*(b)$ and $0\ne c\leq b$ such that $c,h(c),h^2(c)$ are pairwise disjoint.

We show that $f^h$ and $g$ do not commute. Since $f(c)\meet b=0$, $h$ fixes $f(c)$. Also $g(h(c))\meet b=0$, so $h^{-1}$ fixes $g(h(c))$. Using $fg=gf$, we obtain
\[
 gf^h(h(c))=gf(c)=fg(c),\qquad f^h g(h(c))=hfg(h(c)).
\]
In the first case, $fg(b)\meet b=0$, so $h$ fixes $fg(h(c))$. Consequently the two expressions above are $fg(c)$ and $fg(h(c))$. Equality would imply $c=h(c)$, contrary to disjointness. In the second case, $fg$ fixes every subelement of $b$, so the two expressions are $c$ and $h^2(c)$. These are disjoint and nonzero, again a contradiction. This proves (c).

(d) Choose $f\in G$ with $\var(f)=a$, which is possible by Proposition~\ref{prop:support-properties}(b) (or take the identity when $a=0$, a case excluded here). Apply (c) to $f$ and $g$ to obtain $h\in\Sp_G(a)$ such that $f^h$ does not commute with $g$. Since $h(a)=a$, equation \eqref{eq:conj-support} gives $\var(f^h)=a$, so $k=f^h$ has the required properties.
\end{proof}

\begin{proof}[Proof of Theorem~\ref{thm:support-order}]
We first prove
\begin{equation}\label{eq:V-support}
 V(f)=\Sp_G(\var(f)).
\end{equation}
Let $u\in\Sp_G(\var(f))$ and $f'\in D_1(f)$. By Lemma~\ref{lem:phi1-characterization}(b), $\phi_1(f,f')$ implies
\[
 \var(f)\meet\var((f')^{12})=0.
\]
Thus $u$ and $(f')^{12}$ have disjoint variations and commute. This holds for every $f'\in D_1(f)$, so $u\in V(f)$. Hence $\Sp_G(\var(f))\subseteq V(f)$.

Conversely, suppose $u\notin\Sp_G(\var(f))$. Then
\[
 d=\var(u)\meet\neg\var(f)\ne0.
\]
Choose $f'\in\Sp_G^*(d)$ with $(f')^{12}\ne\id$, using Lemma~\ref{lem:local-moving-tools}(b), and put $c=\var((f')^{12})\ne0$. Since $c\leq d\leq\var(u)$, Lemma~\ref{lem:local-moving-tools}(c), applied to $(f')^{12}$ and $u$, gives $h\in\Sp_G(c)$ such that $((f')^{12})^h$ does not commute with $u$. Because $\var(h)\leq c\leq d$, the automorphism $h$ preserves $d$, so for $x=(f')^h$ we have $\var(x)=h(d)=d$. In particular, $\var(x)\meet\var(f)=0$, and Lemma~\ref{lem:phi1-characterization}(a) gives $x\in D_1(f)$. But
\[
 x^{12}=((f')^{12})^h
\]
does not commute with $u$, whence $u\notin V(f)$. This proves the reverse inclusion in \eqref{eq:V-support}.

If $\var(f)\leq\var(g)$, then \eqref{eq:V-support} gives
\[
 V(f)=\Sp_G(\var(f))\subseteq\Sp_G(\var(g))=V(g),
\]
so $\phi_{\leq}(f,g)$ holds. If $\var(f)\nleq\var(g)$, let $a=\var(f)\meet\neg\var(g)\ne0$ and choose $h\in\Sp_G^*(a)$. Then $h\in V(f)$ but $h\notin V(g)$, so $V(f)\nsubseteq V(g)$. This proves the theorem.
\end{proof}
\section*{Acknowledgments}

The author used ChatGPT in preparing
this manuscript.  The author is responsible for the mathematical content of the manuscript.

\end{document}